\documentclass[a4paper,10pt]{article}           
\date{}%

\usepackage{amsmath,amssymb,latexsym}

\newtheorem{defin}{Definition}[section]
\newtheorem{theo}{Theorem}[section]
\newtheorem{lem}{Lemma}[section]
\newtheorem{prop}{Proposition}[section]
\newtheorem{cor}{Corollary}[section]
\newtheorem{rem}{Remark}[section]

\renewcommand{\l} {\pmb{ \big{[}}}
\renewcommand{\r} {\pmb{ \big{]}}}

\title{On some embedding of recursively presented Lie algebras}

\author{E.Chibrikov }

\begin{document}

\baselineskip=1.2\baselineskip

\maketitle

\section { Introduction}

In 1961 G.Higman \cite{Higman} proved an important Embedding Theorem which states that every recursively presented group can be embedded in a finitely presented group. Recall that a group (or an algebra) is called recursively presented if  it can be given by a finite set of generators and a recursively enumerable set of defining relations. If a group (or an algebra) can be given by finite sets of generators and defining relations it is called finitely presented.
As a corollary to this theorem G.Higman proved the existence of a~universal finitely presented group containing every finitely presented group as a subgroup. In fact, its finitely generated subgroups are exactly the finitely generated recursively presented groups.

In \cite{Belyaev} V.Ya.Belyaev proved an analog of Higman's theorem  for associative algebras over a field which is a~finite extention of its simple subfield. The proof was based on his theorem stating that every recursively presented associative algebra over a field as above can be embedded in a recursively presented associative algebra with defining relations which are equalities of words of generators and  $\alpha+\beta=\gamma$, where $\alpha, \beta, \gamma$ are generators. This result allowed Belyaev to apply V.Murskii analog of Higman's theorem \cite{Murskii} for semigroups to obtain a~Higman embedding for associative algebras. In recent paper \cite{Bahturin} Y.Bahturin and A.Olshanskii showed that such embedding can be performed distortion-free. The idea of transition from  algebras to semigroups was also used by G.P.Kukin in \cite {Kukin} (see also \cite{Bokut}). 

This paper appears as a byproduct of the author's joint attemps with Prof. Y.Bahturin to prove Lie algebra analog of Higman's Theorem. In particular, Y.Bahturin suggested to prove that any recursively presented Lie algebra can be embedded in a Lie algebra given by Lie relations of the type mentioned above. In this paper we show that this is indeed true. Namely, every recursively presented Lie algebra over a field which is a~finite extention of its simple subfield can be embedded in a recursively presented Lie algebra defined by relations which are equalities of (nonassociative) words of generators and  $\alpha+\beta=\gamma$ ($\alpha, \beta, \gamma$ are generators). To prove this embedding we use Grobner-Shirshov basis theory for Lie algebras. A short review of this theory as well as some properties of Lyndon-Shirshov words including a new proof of Kukin's Lemma (see \cite{Bokut}, Lemma 2.11.15) are given in the section 2. Note that an existence of Higman's embedding for Lie algebras is still an open problem (see \cite{Sapir}). It is worth to mention a result by L.A.Bokut \cite{Bokut1} that 
for every recursively enumerable set $M$ of positive integers, the Lie algebra 
$$ L_M=\text{Lie}<a,b,c \, | \, [ab^nc]=0, n \in M>,$$
where $[ab^nc]=[[\ldots [[ab]b] \ldots b]c]$, can be embedded into a finitely presented Lie algebra.

This paper was written when the author was a postdoctoral fellow at Memorial University of Newfoundland, supported by a PDF grant by Atlantic Association for Research in Mathematical Sciences as well as Discovery grants of Drs. Y.Bahturin, M.Kotchetov and M.Parmenter. It gives me a pleasure to thank them for the support and encouragement.

\section{ Some definitions and results}

Let $X=\{x_{i}|i\in I\}$ be a linearly ordered set, ${\bf k}$ be a field and
$Lie_{\bf k}(X)$ be the free Lie algebra over ${\bf k}$ generated
by $X$. Define  $\langle X\rangle$ to be the free monoid of all associative
words in $X$ (including the empty word 1).
We use two linear orderings of $\langle X \rangle$:

(i) ({\sl lexicographical order}) $u<1$ for every nonempty word
$u$, and, by induction, $u<v$ if $u=x_{i}u^{\prime}$,

$v=x_{j}v^{\prime}$ and either
$x_{i}<x_{j}$ or $x_{i}=x_{j}$ and $u^{\prime}<v^{\prime}$;

(ii) ({\sl deg-lex order}) $u \prec v$ if $|u|<|v|$, or $|u|=|v|$ and
$u<v$, where $|w|$ denotes the length of $w$.

We regard $Lie_{\bf k}(X)$ as the subspace of the free associative algebra
${\bf k} \langle X \rangle$ which is generated as an~algebra by $X$ under the Lie
bracketing $[xy]=xy-yx$. Given $f \in {\bf k} \langle X \rangle$, denote
by $\bar{f}$ the maximal associative word of $f$ under the ordering (ii);
$f$ is a {\sl monic} if $f=\bar{f}+\sum \alpha_{i}v_{i}$, where
$\alpha_{i} \in {\bf k}$, $v_{i} \in \langle X \rangle$ and
$v_{i} \prec \bar{f}$.


\begin{defin}(\cite{Shirshov}, see also \cite{Shirshov2}, \cite{Lyndon1})
An associative word $w$ is an associative Lyndon--Shirshov
 word if, for arbitrary nonempty $u$ and $v$ such that $w=uv$,  we have $w>vu$.
\end{defin}

We will use the following properties of associative Lyndon-Shirshov words (see \cite{Shirshov}, \cite{Shirshov2} \cite{BokutKolesnikov}, \cite{Bokut}, \cite{Reutenauer}): \\

({\bf 1}) if $w$ is an associative Lyndon-Shirshov word,
then an arbitrary proper prefix of $w$ cannot be 

a~suffix of $w$;

({\bf 2}) if $w=uv$ is an associative Lyndon-Shirshov word, where $u,v \ne 1$, then $w>v$; 

({\bf 3}) if $u$, $v$ are associative Lyndon-Shirshov words and $u>v$, then $uv$ is also an associative 

Lyndon-Shirshov word; 

({\bf 4}) an arbitrary associative word $w$ can be uniquely represented as
$$w=c_{1}c_{2} \ldots c_{n},$$

where $c_{1}, \ldots ,c_{n}$ are  associative Lyndon-Shirshov words and $c_{1} \leq c_{2} \leq \ldots \leq c_{n}$. \\

\begin{defin}
A nonassociative word $[u]$ is a Lyndon-Shirshov word if

(1) $u$ is an associative Lyndon-Shirshov word;

(2) if $[u]=[[u_{1}][u_{2}]]$, then $[u_{1}]$ and $[u_{2}]$ are
Lyndon-Shirshov words (from (1) it then follows that
$u_{1}>u_{2}$);

(3) if $[u]=[[[u_{11}][u_{12}]][u_{2}]]$, then $u_{12} \leq u_{2}$.

\end{defin}

Put  $[x]=x$ for $x \in X$ and  put by induction $[w]=[[u][v]]$ for an associative Lyndon-Shirshov word $w$, where $v$ is the longest proper associative Lyndon-Shirshov end of $w$ (then $u$ is also an associative Lyndon-Shirshov word). Then $[w]$ is a (nonassociative) Lyndon-Shirshov word. The main property of $[w]$ is $$\overline{[w]}=w.$$

It was shown in \cite{Shirshov}, \cite{Chen} (see also \cite{Shirshov2} \cite{Bokut},  \cite{Reutenauer}) that the
set of all  Lyndon-Shirshov words in  the~alphabet $X$ forms a linear basis of $Lie_{\bf k}(X)$.
This implies  that  if $f \in Lie_{\bf k}(X)$ then $\bar{f}$  is an associative Lyndon-Shirshov word. \\

{\bf Lemma} (Shirshov \cite{Shirshov}, see also \cite{Shirshov2})
{\sl Suppose that $w=aub$, where $w$, $u$  are associative Lyndon-Shirshov words. Then a pair of related brackets
in the nonassociative Lyndon-Shirshov word $[w]$ is given by $[uc]$, i.e.
$$
[w]=[a[uc]d],
$$
where $[uc]$ is a nonassociative Lyndon-Shirshov word and $b=cd$.
Represent $c$ in the form
$$
c=c_{1}c_{2} \ldots c_{k},
$$
where $c_{1} \leq c_{2} \leq \ldots \leq c_{k}$ are associative Lyndon-Shirshov words.
Then replacing $[uc]$ by
$[\ldots[[u][c_{1}]] \ldots [c_{k}] ]$
we obtain the word $[w]_{u}$ such that
$$
\overline{[w]}_{u}=w.
$$
}
The word $[w]_{u}$ is called the  special bracketing or the special Shirshov bracketing of $w$ relative to  $u$.\\

We provide some additional information about bracketings in nonassociative Lyndon-Shirshov words.


\begin{lem}
\label{lemma1}
 Take some associative word $w$ represented as
$$
w=c_{1}c_{2} \ldots c_{n},
$$
where $c_{1} \leq c_{2} \leq \ldots \leq c_{n}$ are associative Lyndon-Shirshov words and an associative Lyndon-Shirshov word $v$ which is a~subword of $w$. Then $v$ is a~subword
(possibly not proper) of one of the words $c_{1}$, $c_{2}$,$\ldots$,$c_{n}$.
\end{lem}

{\bf Proof.} Assume that $v$ is not a~subword of any of the words $c_{1}, \ldots, c_{n}$. Then
$$v=c_{i}^{\prime \prime} c_{i+1} \ldots c_{t}^{\prime}$$
for some $1 \leq i < t \leq n$, where $c_{i}=c_{i}^{\prime} c_{i}^{\prime \prime}$, \ $c_{t}=c_{t}^{\prime} c_{t}^{\prime \prime}$ and $c_{i}^{\prime \prime} \ne 1$, $c_{t}^{\prime} \ne 1$. It is possible that $c_{i}=c_{i}^{\prime \prime}$ and/or $c_{t}=c_{t}^{\prime}$. Since a proper prefix is bigger than the word $v$, we have $v < c_{i}^{\prime \prime}$ and $c_{t} \leq c_{t}^{\prime}$.
Property~{\bf 2} of associative Lyndon-Shirshov words implies that
$c_{i}^{\prime \prime} \leq c_{i}$. Since $c_{i} \leq \ldots \leq c_{t}$ it follows that
$v < c_{t}^{\prime}$, but property~2 also yields $v > c_{t}^{\prime}$. 

\vspace{0.3cm}

\begin{cor}
\label{ttt}
Suppose that $w=aub$, where $w$ and $u$  are associative Lyndon-Shirshov words. Then a~pair of related brackets in the nonassociative Lyndon-Shirshov word $[w]$ is given by $[uc]$, i.e.
$$
[w]=[a[uc]d],
$$
where $[uc]$ is a nonassociative Lyndon-Shirshov word and $b=cd$. Represent $c$ in the form
$$
c=c_{1}c_{2} \ldots c_{k},
$$
where $c_{1} \leq c_{2} \leq \ldots \leq c_{k}$ are associative Lyndon-Shirshov words. Then $[uc]=[u[c_{1}][c_{2}] \ldots [c_{k}]]$, that is
$$
[w]=[a[u[c_{1}] \ldots [c_{k}]]d].
$$
\end{cor}

{\bf Proof.}  Shirshov's  Lemma implies that  $[w]=[a[uc]d]$, where $b=cd$ (possibly  $c=1$). Represent the word $c$ as  $c=c_{1}c_{2} \ldots c_{k}$, where $c_{1} \leq c_{2} \leq \ldots \leq c_{k}$ are associative Lyndon-Shirshov words. Let us consider the Lyndon-Shirshov word $[uc]=[uc_{1}c_{2} \ldots c_{k}]$. By applying  Shirshov's Lemma and Lemma \ref{lemma1} to the words $c_{1}$, $c_{2}$,\ldots ,$c_{k}$  we deduce that
 $[uc]=[u[c_{1}][c_{2}] \ldots [c_{k}]]$. \\

 The following Lemma is a generalization of the Shirshov Lemma and for the first time was noticed by G.Kukin (see \cite{Bokut}, Lemma 2.11.15). To prove this Lemma we use a different approach.

\begin{lem} (Kukin)
\label{lemma2}
Suppose that $w=aubvd$, where $w$, $u$, $v$ are associative Lyndon-Shirshov words. Then there is some bracketing
$$
[w]_{u,v}=[a[u]b[v]d]
$$
in the word $w$ such that
$$
\overline{[w]}_{u,v}=w.
$$
\end{lem}

{\bf Proof.}  Shirshov's Lemma implies that $[w]=[a[uc]p]$, where $[uc]$ is a nonassociative  Lyndon-Shirshov word and $cp=bvd$.
If $v$ is a subword of $p$, then $[w]=[a[uc]q[vs]l]$ by  the Shirshov  Lemma, where $[vs]$ is a Lyndon-Shirshov word and $p=qvsl$.
Write $c=c_{1}c_{2} \ldots c_{n}$ and $s=s_{1}s_{2} \ldots s_{m}$, where $c_{1} \leq c_{2} \leq \ldots \leq c_{n}$ and  $s_{1} \leq s_{2} \leq \ldots \leq s_{m}$ are associative Lyndon-Shirshov words.
 Then we define
$$
[w]_{u,v}=[a \ [[u][c_{1}] \ldots [c_{n}] ] \ q   \ [[v][s_{1}]  \ldots [s_{m}] ] \ l],
$$
where
$$[[u][c_{1}] \ldots [c_{n}] ]=[ \ldots [[u][c_{1}]] \ldots [c_{n}] ]
$$
and
$$[[v][s_{1}] \ldots [s_{m}] ]=[ \ldots [[v][s_{1}]] \ldots [s_{m}] ]
$$
are Shirshov's special bracketings  in the words $uc$ and $vs$ relative to $u$ and $v$ respectively (see Shirshov's Lemma). Therefore, we get
$$
\overline{[w]}_{u,v}=w.
$$

Assume now that the word $v$  is not a subword of $p$. Shirshov's Lemma implies that $[w]=[aub[vs]q]$, where $[vs]$ is a nonassociative Lyndon-Shirshov word and $sq=d$. Therefore, $[vs]$ is a nonassociative subword of $[uc]$, hence $v$ is a subword of $c$.

Represent the word $c$ as  $c=c_{1}c_{2} \ldots c_{n}$, where  $c_{1} \leq c_{2} \leq \ldots \leq c_{n}$ are associative Lyndon-Shirshov words. By Lemma \ref{lemma1} we deduce that  $v$ is a subword of $c_{t}$ for some $t$.  If $v=c_{t}$, then we define
$$
[w]_{u,v}=[w]_{u}= [a \  [  [u]  \ [c_{1}]  \ldots  [c_{t}]  \ldots [c_{n}]  ] \  p ],
$$
i.e. $[w]_{u,v}$ is the special bracketing in $w$ relative to $u$. It is obvious that
$$
\overline{[w]}_{u,v}=w.
$$
Now, we suppose that $v$ is a proper subword of $c_{t}$. Corollary \ref{ttt} gives us the equality:
$$
[w]=[a \ [ u \ [c_{1}] \ldots [c_{t}] \ldots [c_{n}] \ ] \ p ].
$$
Consequently, $vs$ is a subword of $c_{t}$, that is $c_{t}=c_{t}^{\prime} vs c_{t}^{\prime \prime}$ for some $c_{t}^{\prime}$ and $c_{t}^{\prime \prime}$.
Write  $s=s_{1} s_{2} \ldots s_{m}$, where $s_{1} \leq s_{2} \leq \ldots \leq s_{m}$ are associative Lyndon-Shirshov words. Then we define
$$
[w]_{u,v}=[a \ [ [u] \ [c_{1}] \ldots [c_{t}]_{v} \ldots [c_{n}]  ] \ p ],
$$
where $[c_{t}]_{v}=[c_{t}^{\prime} \  [[v] [s_{1}] \ldots [s_{m}]  ] \ c_{t}^{\prime \prime} ]$ is the  Shirshov special bracketing in $c_{t}$ relative to $v$.  The word $[w]_{u,v}$ results from  $[ [u] \ [c_{1}] \ldots [c_{t}] \ldots [c_{n}]  ]$, which is the Shirshov special bracketing in $uc_{1} \ldots c_{t} \ldots c_{n}$ relative to $u$, by replacing $[c_{t}]$ by $[c_{t}]_{v}$. It is clear that
$$
\overline{[w]}_{u,v}=w.
$$

\begin{rem}
It is not difficult to see that Lemma \ref{lemma2} can be easily generalized to an arbitrary number of associative Lyndon-Shirshov subwords of $w$.
\end{rem}

\vspace{0.3cm}


Let $f$ be a monic Lie polynomial such that $\bar{f}=h$ is a subword of $w$, i.e. $w=u \bar{f} v$.
Denote by $[ufv]_{f}$  the Lie polynomial obtained from $[uhv]_{h}$ by replacing $[h]$ by $f$. The polynomial $[ufv]_{f}$ has the property
 $$\overline{[ufv]}_{f}=w.$$

\begin{defin}
Take monic Lie polynomials $f$ and $g$ and a word $w$ such that $w=\bar{f}u=v\bar{g}$, where $u,v \in \langle X \rangle$ and $|\bar{f}|+|\bar{g}|>|w|$. The intersection composition of $f$ and $g$ relative
to $w$ is defined by
$$
(f,g)_{w}=[fu]_{f}-[vg]_{g}.
$$
\end{defin}

\begin{defin}
Take monic Lie polynomials $f$ and $g$ and a word $w$ such that $w=\bar{f}=u \bar{g} v$, where $u,v \in \langle X \rangle$. The inclusion composition of $f$ and $g$ relative to $w$ is defined by
$$
(f,g)_{w}=f-[ugv]_{g}.
$$
\end{defin}

The  main properties of  compositions are:
$(f,g)_{w} \in Id(f,g)$ and $\overline{(f,g)}_{w} \prec w.$

\vspace{0.2cm}

\begin{defin} (\cite{Bokut1})
Given a set $S$ of monic Lie polynomials, the composition $(f,g)_{w}$ of $f$ and $g$ is called trivial relative to $S$ if $(f,g)_{w}=\sum_{i}\alpha_{i}[u_{i}s_{i}v_{i}]_{s_{i}}$, where
$\alpha_{i} \in {\bf k}$, $u_{i}, v_{i} \in \langle X \rangle$,
$s_{i} \in S$ and $u_{i}\bar{s}_{i}v_{i} \prec w$.
\end{defin}

\vspace{0.1cm}

\begin{defin}
Let $S$ be a set of monic Lie  polynomials. Then $S$ is a Gr\"{o}bner-Shirshov basis if every composition of any two elements of $S$ is trivial relative to $S$.
\end{defin}

\vspace{0.1cm}

\begin{rem}
\label{remark1}
If $f \in Lie_{\bf k}(X)$, then $\bar{f}$  is an associative Lyndon-Shirshov word and by  property {\bf 1}
of associative Lyndon-Shirshov words, there is no composition $(f,f)_{w}$ for any $w \in \langle X\rangle$. Therefore, $S=\{ f \}$ is a Gr\"{o}bner-Shirshov basis for every $f \in Lie_{\bf k}(X)$.
\end{rem}

\vspace{0.1cm}

A Lyndon-Shirshov word $[w]$ is called $S$-reduced if $w \ne u \bar{s} v$ for any $s \in S$ and $u,v \in \langle X \rangle$. A set $Y$ is called well-ordered if any chain of elements $w_1> \ldots > w_k > \ldots$ of $\langle Y \rangle$ is stabilized. An~important result is Shirshov's Composition Lemma.
We recall it as stated in \cite{Bokut1}: \\

{\bf Composition Lemma }
{\sl Let $X$ be a well-ordered set.  If $S$ is a Gr\"{o}bner-Shirshov basis and
$f \in Id(S)$, then $\bar{f}=u\bar{s}v$ for some $s \in S$ and $u, v \in \langle X \rangle$.} \\

It has an important corollary: \\

{\bf Composition-Diamond Lemma}
{\sl Suppose that $X$ is a well-ordered set. Then $S$ is a Gr\"{o}bner-Shirshov basis
if and only if the set of all Lyndon-Shirshov $S$-reduced words is a linear basis for $Lie_{\bf k}(X)/ Id(S)=Lie_{\bf k}(X|S)$. } \\


Let $f$ and $g$ be  elements of $Lie_{\bf k}(X)$. In some applications the following question appears:  \\ which conditions on $f$ and $g$
give rise to the equality
\begin{equation}
\label{mmm1}
Id(f) \bigcap Id(g)=Id(f)\cdot Id(g),
\end{equation}
where $Id(h)$ denotes the ideal of $Lie_{\bf k}(X)$ generated by $h$? 

As a corollary to Kukin's Lemma we obtain the following sufficient condition for the equality (\ref{mmm1}) to hold.

\begin{prop}
If there is a linear order of $X$ such that there are no compositions $(f,g)_w$ and $(g,f)_w$ for any $w \in \langle X \rangle$, then  $Id(f) \bigcap Id(g)=Id(f) \cdot Id(g)$.
\end{prop}

{\bf Proof}.   If $h \in Id(f) \bigcap Id(g)$, then
$$
h=\sum_{i} \alpha_{i} [a_{i}fb_{i}]= \sum_{j} \beta_{j}[c_{j}gd_{j}],
$$
where $\alpha_{i}, \beta_{j} \in {\bf k}$ and $[...]$ are some bracketings.
The Composition-Diamond Lemma  (see also Remark~\ref{remark1}) and  the hypotheses of the proposition imply that the maximal
associative word of the polynomial $h$ is
$$
\bar{h}=u_{1}\bar{f}u_{2}\bar{g}u_{3}.
$$
Let us consider
$$
h-[u_{1}f u_{2}g u_{3}]_{f, g},
$$
where $[u_{1} f u_{2} g u_{3}]_{f, g}$ is obtained from $[u_{1}\bar{f}u_{2}\bar{g}u_{3}]_{\bar{f}, \bar{g}}$ by replacing $\bar{[f]}$ and $\bar{[g]}$ by $f$ and $g$ respectively.
Lemma \ref{lemma2} implies that
$$
\overline{[u_{1} f u_{2} g u_{3}]}_{f, g}=\bar{h}.
$$
Since
$$
h-[u_{1} f u_{2} g u_{3}]_{f, g} \in Id(f) \bigcap Id(g)
$$
and
$$
\overline{h-[u_{1} f u_{2} g u_{3}]}_{f, g} < \bar{h},
$$
the claim of the proposition follows by induction. \\


\section{Definition of the Functions $n_i(j)$ and $s(i,j)$}

The idea of these functions belongs to V.Ya.Belyaev. We introduce them slightly modified. 

Suppose that $\phi(i,j)$ is a recursive function. We define two functions 
$$
n_i(j): \  \mathbb{N} \times \mathbb{N} \rightarrow \mathbb{N} \text{ and }
s(i,j): \  \mathbb{N} \times \mathbb{N} \rightarrow \mathbb{N}
$$
as follow: 

Let $i>j$ and $n_2(1)$ be a natural number which is greater than or equal to $max\{ 2, \phi(2,1). \}$. Then we define $n_i(j)$ as an arbitrary natural number such that $n_i(j) \geq max \{ n_i(j-1)+1, \phi(i,j) \}$ for $j \geq 2$ and $n_i(1) \geq max \{ n_{i-1}(i-2)+1, \phi(i,1) \}$ for $j=1$. For $i \leq j$ we set $n_i(j)=n_{j+1}(i)$.

Let $i<j$. We define $s(i,j)=s(j,i)=n_j(i)$.

It is not difficult to notice that the functions satisfy the following properties:
$$
n_{i}(j) \geq i \text{ for } \forall j \text{ and } s(i,j) \geq \phi(i,j).
$$

\section{Encoding of Addition and Multiplication}

Now, let $X=\{ x, y, z \}$, ${\bf k}$ be a field which is a finite extension of its simple subfield and let $Lie_{\bf k}(X)$ be the~free Lie algebra over ${\bf k}$ generated by $X$. If $w=[..[[x_1,x_2],x_3], \ldots, x_n]$ is a left-normed word, we will write $w=\l x_1x_2 \ldots x_n \r_l$. Denote by $I$ the ideal of $Lie_{\bf k}(X)$ generated by the following relations:
\begin{equation}
\label{eq1}
\l xy^iz \r_l + \l xy^jz \r_l- \l xy^{\phi(i,j)}z \r_l, \ i>j; \ i,j=1,2, \ldots
\end{equation}
\begin{equation}
\label{eq2}
\big{[} \l xy^iz \r_l, \l xy^jz \r_l \big{]}  - \l xy^{\psi(i,j)}z \r_l, \ i>j; \ i,j=1,2, \ldots
\end{equation}
\begin{equation}
\label{eq3}
\delta_t \l xy^iz \r_l - \l xy^jz \r_l, \ (i,j) \in D_1,
\end{equation}
\begin{equation}
\label{eq4}
\l xy^iz \r_l - \l xy^jz \r_l, \ (i,j) \in D_2,
\end{equation} 
where $\phi$, $\psi$ are recursive functions, $D_1$, $D_2$ are recursively denumerable sets, $\delta_t$ $(t=1,2,\ldots,n)$ are generators of the field ${\bf k}$. 

Let $X_1=X \bigcup \{ u, \alpha, \beta, \gamma \}$ and $Lie_{\bf k} (X_1)$ be the free Lie algebra over the field ${\bf k}$ generated by $X_1$. Denote by $I_1$ the ideal of $Lie_{\bf k} (X_1)$ generated by the relations:

\begin{equation}
\label{eq5}
\l xy^iz \r_l - \l xu^{n_i(j)} \epsilon_{ij} z \r_l; \ i,j=1,2, \ldots; \ 
\epsilon_{ij}=
\left\{
\begin{array}{ll}
\alpha, & i+j \text{ is even },  \\
\beta,  & i+j \text{ is odd },
\end{array} \right.
\end{equation}
\begin{equation}
\label{eq6}
\l xy^{\phi(i,j)}z \r_l - \l xu^{s(i,j)} \gamma z \r_l, \ i \ne j; \ i,j=1,2, \ldots
\end{equation}
\begin{equation}
\label{eq7}
\alpha+\beta-\gamma,
\end{equation}
\begin{equation}
\label{eq8}
\big{[} \l xy^iz \r_l, \l xy^jz \r_l \big{]}  - \l xy^{\psi(i,j)}z \r_l, \ i>j; \ i,j=1,2, \ldots
\end{equation}
\begin{equation}
\label{eq9}
\delta_t \l xy^iz \r_l - \l xy^jz \r_l, \ (i,j) \in D_1,
\end{equation}
\begin{equation}
\label{eq10}
\l xy^iz \r_l - \l xy^jz \r_l, \ (i,j) \in D_2,
\end{equation} 
where $\phi$, $\psi$ are the recursive functions, $D_1$, $D_2$ are the recursively denumerable sets, $\delta_t$ $(t=1,2,\ldots,n)$ are the generators of the field ${\bf k}$ and $n_i(j)$, $s(i,j)$ are functions defined in 1.1 relative to $\phi$.

\begin{lem}
\label{L1}
$I \subseteq I_1$.
\end{lem}

{\bf Proof}. It suffices to show that the relation (\ref{eq1}) belogns to $I_1$. Let $i<j$, then in $Lie_{\bf k}(X_1)/I_1$ we have 
$$
\l xy^iz \r_l+ \l xy^jz \r_l=\l x u^{n_i(j-1)} \epsilon_{i(j-1)} z \r_l +
\l x u^{n_j(i)} \epsilon_{ij}z \r_l. 
$$
Since $n_i(j-1)=n_j(i)$ we obtain that
$$
\l xy^iz \r_l+ \l xy^jz \r_l=\l x u^{s(i,j)}(\epsilon_{i(j-1)}+\epsilon_{ij})z \r_l=
\l x u^{s(i,j)}(\alpha + \beta)z \r_l=\l x u^{s(i,j)} \gamma z \r_l= 
\l x y^{\phi(i,j)}z \r_l.
$$
Thus, the relation $\l xy^iz \r_l+ \l xy^jz \r_l - \l x y^{\phi(i,j)}z \r_l$ belogns to $I_1$. \\

Let $x > z > \alpha > \beta > \gamma > u > y$. Denote by $R$ a set of relations the union of which and (\ref{eq1})-(\ref{eq4}) is closed with respect to the operation of composition, that is the union is a Grobner-Shirshov basis.
It follows from the definition of composition that $R \subseteq I$ and the set $R$ consists of Lie polynomials in $\l xy^iz \r_l$. \\

\begin{lem}
The following set $S$ forms a Grobner-Shirshov basis:
\begin{equation}
\label{eq11}
\l xu^{n_i(j)} \epsilon_{ij} z \r_l - \l xy^iz \r_l; \ i,j=1,2, \ldots; \ 
\epsilon_{ij}=
\left\{
\begin{array}{ll}
\alpha, & i+j \text{ is even },  \\
\beta,  & i+j \text{ is odd },
\end{array} \right.
\end{equation}
\begin{equation}
\label{eq12}
\l xu^{s(i,j)} \gamma z \r_l - \l xy^{\phi(i,j)}z \r_l, \ i \ne j; \ i,j=1,2, \ldots
\end{equation}
\begin{equation}
\label{eq13}
\alpha+\beta-\gamma,
\end{equation}
\begin{equation}
\label{eq135}
\l xy^iz \r_l + \l xy^jz \r_l- \l xy^{\phi(i,j)}z \r_l, \ i>j; \ i,j=1,2, \ldots
\end{equation}
\begin{equation}
\label{eq14}
\big{[} \l xy^iz \r_l, \l xy^jz \r_l \big{]}  - \l xy^{\psi(i,j)}z \r_l, \ i<j; \ i,j=1,2, \ldots
\end{equation}
\begin{equation}
\label{eq15}
\delta_t \l xy^iz \r_l - \l xy^jz \r_l, \ (i,j) \in D_1,
\end{equation}
\begin{equation}
\label{eq16}
\l xy^iz \r_l - \l xy^jz \r_l, \ (i,j) \in D_2,
\end{equation} 
\begin{equation}
\label{eq17}
R
\end{equation}
where $\phi$, $\psi$ are the recursive functions, $D_1$, $D_2$ are the recursively denumerable sets, $\delta_t$ $(t=1,2,\ldots,n)$ are the generators of the field ${\bf k}$.
\end{lem}

{\bf Proof}. Note that all the words of the relations (\ref{eq11})-(\ref{eq17}) are nonassociative Lyndon-Shirshov words. Since $n_i(j) \geq i$ for $\forall j$ and $s(i,j) \geq \phi(i,j)$, the leading words of the relations (\ref{eq11})-(\ref{eq13}) contain  $\alpha$, $\beta$ or $\gamma$. Therefore, the relations (\ref{eq11})-(\ref{eq13}) do not form compositions with the relations (\ref{eq135})-(\ref{eq17}). Because the~relations (\ref{eq135})-(\ref{eq17}) are closed under composition by the definition of $R$, we have only one composition to verify that it is trivial, namely, the inclusion composition of  (\ref{eq11}) and (\ref{eq13}).

Let $f= \l x u^{n_i(j)} \alpha z \r_l - \l x y^i z \r_l, \text{ where } i+j \text{ is even }$, and $ g=\alpha + \beta - \gamma$. Then $w=x u^{n_i(j)} \alpha z$ and
$$
(f,g)_w = f - \l x u^{n_i(j)} ( \alpha + \beta - \gamma) z \r_l = - \l x y^i z \r_l - 
\l x u^{n_i(j)} \beta z \r_l + \l x u^{n_i(j)} \gamma z \r_l.
$$
If $i \leq j$, then $n_i(j)=n_{j+1}(i)$ and
$$
(f,g)_w = -( \l x u^{n_{j+1}(i)} \beta z \r_l - \l x y^{j+1} z \r_l ) + 
( \l xu^{s(i,j+1)} \gamma z \r_l - \l x y^{\phi(i,j+1)} z \r_l )- ( \l x y^i z \r_l +  \l x y^{j+1} z \r_l - \l x y^{\phi(i,j+1)} z \r_l ). 
$$
If $i>j$, then $n_j(i-1)=n_i(j)$ and 
$$
(f,g)_w = -( \l x u^{n_j(i-1)} \beta z \r_l - \l x y^j z \r_l ) + 
( \l x u^{s(i,j)} \gamma z \r_l - \l x y^{\phi(i,j)} z \r_l )- ( \l x y^i z \r_l +  \l x y^j z \r_l - \l x y^{\phi(i,j)} z \r_l ). 
$$
Therefore, $(f,g)_w$ is trivial relative to $S$. \\

\begin{prop}
\label{prop1}
$Lie_{\bf k}(X)/I$ is embedded into $Lie_{\bf k}(X_1)/I_1 $.
\end{prop}

{\bf Proof}. It is sufficient to prove that $Lie_{\bf k}(X) \bigcap I_1=I$. The inclusion $I \subseteq Lie_{\bf k}(X) \bigcap I_1$ follows from Lemma \ref{L1}. Suppose that $f \in Lie_{\bf k}(X) \bigcap I_1$. Composition Lemma for Lie algebras implies that in the~free associative algebra ${\bf k} \langle X_1 \rangle$  \  \   $\overline{f}=a \overline{s} b$ for some words $a,b \in \langle X_{1} \rangle$ and $s \in S$ (recall that $I_1=Id(S)$). But $f \in Lie_{\bf k} (X)$, hence, $\overline{f} \in \langle X \rangle$ and, in particular, $\overline{s} \in \langle X \rangle$. Since the relations (\ref{eq11})-(\ref{eq13}) contain $\alpha$, $\beta$ or $\gamma$ in their leading words, $s$ is a relation of (\ref{eq135})-(\ref{eq17}) and, therefore, $s \in I$. We have $f - [ a s b ]_{\overline{s}} \in Lie_{\bf k} (X) \bigcap I_1$ and $\overline{f - [ a s b ]_{\overline{s}}} < \overline{f}$. By performing the same calculations at some step we obtain that
$$
f - [asb]_{\overline{s}}-[a_1s_1b_1]_{\overline{s_1}}-...-[a_ks_kb_k]_{\overline{s_k}}=0,
$$
where $s, s_1,..., s_k \in I$ and $a,b,a_1,b_1,...,a_k,b_k \in \langle X \rangle$. Thus, $f \in I$.

\section{Proof of Embedding}

Let $Y=\{ a, b, c \}$ with $a>c>b$ and $Lie_{\bf k}(Y)$ be the free Lie algebra generated by $Y$ over a field ${\bf k}$ which is a finite extention of its simple subfield. Note that the set of left-normed elements 
$$
T=\{ \l a b ^i c \r_l \ | \ i=1,2,3,\ldots \}
$$ 
generates a free Lie subalgebra and constitutes a set of free generators. Indeed, the words $\l a b ^i c \r_l$ are Lyndon-Shirshov words and by ordering $T$ as $\l a b ^i c \r_l > \l a b ^j c \r_l$ if and only if $i<j$, we obtain that every nonassocitive Lyndon-Shirshov word in the alphabet $T$ is also a Lyndon-Shirshov word in the alphabet $Y$. It implies that Lyndon-Shirshov words in the alphabet $T$ are linearly independent in $Lie_{\bf k}(Y)$, thus, $T$ is a set of free generators.

Let $L$ be a finitely generated Lie algebra over the field ${\bf k}$ with a recursively denumerable set of defining relations. Then $L$ can be represented as  $L=\{ a_1, a_2, a_3, \ldots \}$, where each element of $L$ is repeated at least twice. There exists an embedding $\Phi: L \longrightarrow Lie_{\bf k}(Y)/J$ such that  $\Phi(a_i)= \l a b^i c \r_l$ and the ideal $J$ of $Lie_{\bf k}(Y)$ is generated by the following relations:
\begin{equation}
\label{eq19}
\l ab^ic \r_l + \l ab^jc \r_l- \l ab^{\phi(i,j)}c \r_l, \ i>j; \ i,j=1,2, \ldots
\end{equation}
\begin{equation}
\label{eq20}
\big{[} \l ab^ic \r_l, \l ab^jc \r_l \big{]}  - \l ab^{\psi(i,j)}c \r_l, \ i>j; \ i,j=1,2, \ldots
\end{equation}
\begin{equation}
\label{eq21}
\delta_t \l ab^ic \r_l - \l ab^jc \r_l, \ (i,j) \in D_1,
\end{equation}
\begin{equation}
\label{eq22}
\l ab^ic \r_l - \l ab^jc \r_l, \ (i,j) \in D_2,
\end{equation} 
where $\phi$, $\psi$ are recursive functions, $D_1$, $D_2$ are recursively denumerable sets, $\delta_t$ $(t=1,2,\ldots,n)$ are generators of the field ${\bf k}$. 

Indeed, the Lie algebra $Lie_{\bf k}(T)/J^{\prime}$, where $J^{\prime}$ is an ideal of $Lie_{\bf k}(T)$ generated by (\ref{eq19})-(\ref{eq22}), is isomorphic to $L$. We have to show that $Lie_{\bf k}(T) \bigcap J=J^{\prime}$. The inclusion $J^{\prime} \subseteq Lie_{\bf k}(T) \bigcap J$ is obvious. Let 
$f \in Lie_{\bf k}(T) \bigcap J$. By Composition Lemma we have $\overline{f}=a \overline{s} b$, where $a,b \in \langle Y \rangle$ and $s$ is a composition obtained from the relations (\ref{eq19})-(\ref{eq22}). Since $s$ is a Lie polynomial in $\l ab^i c \r_l$ (hence, $s \in J^{\prime}$) and $f \in {\bf k} \langle T \rangle$ we can write
$$
\overline{f}=ab^{i_1} c \ldots ab^{i_t} c \  \overline{s} \ ab^{i_{t+1}} c \ldots ab^{i_k} c.
$$
Therefore, 
$$
f-[ \l ab^{i_1} c \r_l \ldots \l ab^{i_t} c \r_l \  s \  \l ab^{i_{t+1}} c \r_l \ldots \l ab^{i_k} c \r_l  ]_{\overline{s}} \in Lie_{\bf k}(T) \bigcap J \text{ and }
$$
$$
\overline{ f-[ \l ab^{i_1} c \r_l \ldots \l ab^{i_t} c \r_l \  s \  \l ab^{i_{t+1}} c \r_l \ldots \l ab^{i_k} c \r_l  ]_{\overline{s}} } < \overline{f},
$$
and the proof follows from the induction on the degree of $f$. \\

Now, Proposition \ref{prop1} implies the following

\begin{theo}
A recursively presented Lie algebra over a field which is a finite extention of its simple subfield can be embedded into a recursively presented Lie algebra defined by relations which are equalities of (nonassociative) words of generators and  $\alpha+\beta-\gamma$, where $\alpha, \beta, \gamma$ are generators.
\end{theo}


\begin{thebibliography}{11}

\bibitem{Bahturin} \emph{Y.Bahturin, A.Olshanskii} Filtrations and distortions in infinite-dimensional algebras, J.Algebra, to~appear; doi:10.1016/j.jalgebra.2010.09.019.

\bibitem{Belyaev} \emph{V.Ya.Belyaev} Subrings of finitely presented associative rings, Algera i Logica, {\bf 17} (1978), 627-638.

\bibitem{Bokut1}  \emph{L.A.Bokut } Unsolvability of the equality problem and subalgebras of finitely presented Lie algebras. Math. USSR Izvestia
{\bf 6} (1972), 1153--1199.


\bibitem{BokutKolesnikov} \emph{L.A.Bokut, P.S.Kolesnikov} Gr\"{o}bner-Shirshov bases: from their incipiency to the present. Journal of Mathematical Sciences,
{\bf 116} (2003), 1,  2894--2916.

\bibitem{Bokut} \emph{L.A.Bokut, G.P.Kukin }  Algorithmic and combinatorial
algebra.  Mathematics  and its Applications,  Kluwer Academic Publishers
 Group, Dordrecht, (1994). 


\bibitem{Chen} \emph{K.T.Chen, R.H.Fox, R.C.Lyndon} Free differential
calculus.IV: The  quotient groups of the lower central series. Ann. Math.
{\bf 68} (1958), 2, 81--95.

\bibitem{Higman} \emph{G.Higman} Subgroups of finitely presented groups, Proc. Royal Soc. London (Series A), 262 (1961), 455-475.

\bibitem{Kukin} \emph{G.P.Kukin} On the equality problem for Lie algebras, Sib. Mat. Zh., {\bf 18}, No. 5, 1194-1197 (1977).


\bibitem{Sapir} \emph{O.~G.~Kharlampovich, M.~V.~Sapir} Algorithmic problems in varieties, Inter. J.Algebra Compt., 5, no.4, 5 (1995) 379--602.

\bibitem{Murskii} \emph{V.L.Murskii} Isomorphic embeddability of a semigroup with an enumerable set of defining relations into a finitely presented semigroup, Mat.Zemetki, {\bf 1}(2) (1967), 217-224.

\bibitem{Lyndon1} \emph{R.C.Lyndon } On Burnside's problem I. Trans. Amer.
 Math. Soc. {\bf 77} (1954), 202--215.
 
 \bibitem{Reutenauer} \emph{C.Reutenauer} Free Lie algebras. London
Mathematical Society Monographs. New Series, 7. Oxford Science Publications.
The Clarendon Press, Oxford University Press, New York, (1993).
 
 
\bibitem{Shirshov} \emph{A.I.Shirshov } On free Lie rings.  Mat. Sb.,
{\bf 45} (1958), 2, 113--122.

\bibitem{Shirshov1} \emph{A.I.Shirshov} Sertain algorithmic problems for Lie algebras. Sibirsk. Mat. Z. {\bf 3} (1962), 292--296 (Translation in ACM SIGSAM Bull
{\bf 33} (1999), 2, 3--6).

\bibitem{Shirshov2} \emph{Selected works of A.I.Shirshov}, Birkh\"{a}user, 2009.

\end{thebibliography}
\end{document}